\ifx \labelsONmargin\undefined

  \ifx\RequirePackage\undefined
    \expandafter\ifx\csname amsart.sty\endcsname\relax
        \documentstyle[12pt,amssymb,amscd]{amsart}
    \fi

    \def\atSign{@@}

    \def\mathbb{\Bbb}
    \def\mathfrak{\frak}
    \def\mathbf{\bold}
    \ifx\boldsymbol\undefined   
      \def\boldsymbol#1{{\bold #1}}
    \fi

    \def\mathbit{\boldsymbol}

    \makeatletter
    \newenvironment{proof}{%
         \@ifnextchar[{%
                       \expandafter\let\expandafter\end@proof
                         \csname endpf*\endcsname
                         \my@proof
                      }{\let\end@proof\endpf\pf}%
        }{\end@proof}
    \def\my@proof[#1]{\@nameuse{pf*}{#1}}
    \makeatother

    \def\xrightarrow[#1]#2{@>{#2}>{#1}>}
    \def\xleftarrow[#1]#2{@<{#2}<{#1}<}

    \def\providecommand#1{\def#1}
    \def\emph#1{{\em #1}}
    \def\textbf#1{{\bf #1}}

    \def\mathring{\overset{\,\,{}_\circ}}
  \else
      \ifx\usepackage\RequirePackage
      \else

        \documentclass[12pt]{amsart}
        \usepackage{amssymb,amscd,diagrams}
        
      \fi
      \ifx \mathring\undefined      
        \DeclareMathAccent{\mathring}{\mathalpha}{operators}{"17}
      \fi

      \long\def\FAKEendPROOF{\endtrivlist}
      \ifx\endproof\FAKEendPROOF
      \def\endproof{\qed\endtrivlist}
      \fi

      \ifx\mathbit\undefined    
        \DeclareMathAlphabet{\mathbit}{OML}{cmm}{b}{it}
      \fi

      \def\atSign{@}

      \advance\oddsidemargin -0.1\textwidth
      \evensidemargin\oddsidemargin
      \def\Sb#1\endSb{_{\substack{#1}}}
      \def\Sp#1\endSp{^{\substack{#1}}}

  \fi
\makeatletter
\@ifundefined{DeclareFontShape}
     {\@ifundefined{selectfont}
             {
                \message{We need AmS-LaTeX, this version of LaTeX does not
                        support it}\@@end
             }{
                \def\mathcal{\cal}
                
                \def\pcyr{%
                        \def\default@family{UWCyr}%
                        \let\oldSl@\sl
                        \def\sl{\def\default@shape{it}\oldSl@}%
                        \cyracc
                        \language\Russian\family{UWCyr}\selectfont
                }
             }}{
                \DeclareFontEncoding{OT2}{}{\noaccents@}
                \DeclareFontSubstitution{OT2}{cmr}{m}{n}
                \DeclareFontFamily{OT2}{cmr}{\hyphenchar\font45 }
                \DeclareFontShape{OT2}{cmr}{m}{n}{%
                     <5><6><7><8><9><10>gen*wncyr %
                     <10.95><12><14.4><17.28><20.74><24.88> wncyr10 %
                }{}
                \DeclareFontShape{OT2}{cmr}{m}{it}{%
                     <5><6><7><8><9><10> gen * wncyi%
                     <10.95><12><14.4><17.28><20.74><24.88> wncyi10%
                }{}
                \DeclareFontShape{OT2}{cmr}{bx}{n}{%
                     <5><6><7><8><9><10> gen * wncyb%
                     <10.95><12><14.4><17.28><20.74><24.88> wncyb10%
                }{}
                \DeclareFontShape{OT2}{cmr}{m}{sl}{%
                     <-> ssub * cmr/m/it%
                }{}
                \DeclareFontShape{OT2}{cmr}{m}{sc}{%
                     <5><6><7><8><9><10>%
                     <10.95><12><14.4><17.28><20.74><24.88> wncysc10%
                }{}
                \DeclareFontFamily{OT2}{cmss}{\hyphenchar\font45 }
                \DeclareFontShape{OT2}{cmss}{m}{n}{%
                     <8><9><10> gen * wncyss%
                     <10.95><12><14.4><17.28><20.74><24.88> wncyss10%
                }{}
                
                \def\cyrencodingdefault{OT2}
                \def\pcyr{%
                        \cyracc
                        \let\encodingdefault\cyrencodingdefault
                        \language\Russian\fontencoding{OT2}\selectfont
                }
             }

\@ifundefined{theorembodyfont}
     {
        \def\theorembodyfont#1{\relax}
        \makeatletter
          \let\@@th@plain\th@plain
          \def\th@plain{ \@@th@plain \slshape }
        \makeatother
        \let\normalshape\relax
     }{}

\ifx\cprime\undefined
     \def\cprime{$'$}
\fi
\makeatletter
  \def\@sect@my#1#2#3#4#5#6[#7]#8{%
\ifnum #2>\c@secnumdepth
   \let\@svsec\@empty
 \else
   \refstepcounter{#1}%
\edef\@svsec{\ifnum#2<\@m
             \@ifundefined{#1name}{}{\csname #1name\endcsname\ }\fi
\noexpand\rom{\csname the#1\endcsname.}\enspace}\fi
 \@tempskipa #5\relax
 \ifdim \@tempskipa>\z@ 
   \begingroup #6\relax
   \@hangfrom{\hskip #3\relax\@svsec}{\interlinepenalty\@M #8\par}%
   \endgroup
   \if@article\else\csname #1mark\endcsname{%
        \ifnum \c@secnumdepth >#2\relax\csname the#1\endcsname. \fi#7}\fi
\ifnum#2>\@m \else
       \let\@tempf\\ \def\\{\protect\\}\addcontentsline{toc}{#1}%
{\ifnum #2>\c@secnumdepth \else
             \protect\numberline{%
               \ifnum#2<\@m
               \@ifundefined{#1name}{}{\csname #1name\endcsname\ }\fi
               \csname the#1\endcsname.}\fi
           #8}\let\\\@tempf
     \fi
 \else
  \def\@svsechd{#6\hskip #3\@svsec
    \@ifnotempty{#8}{\ignorespaces#8\unskip
       \ifnum\spacefactor<1001.\fi}%
        \ifnum#2>\@m \else
          \let\@tempf\\ \def\\{\protect\\}\addcontentsline{toc}{#1}%
            {\ifnum #2>\c@secnumdepth \else
              \protect\numberline{%
                \ifnum#2<\@m
                \@ifundefined{#1name}{}{\csname #1name\endcsname\ }\fi
                \csname the#1\endcsname.}\fi
             #8}\let\\\@tempf\fi}%
 \fi
\@xsect{#5}}
  \ifx\RequirePackage\undefined
  \let\@sect\@sect@my             
  \fi
  \ifx\RequirePackage\undefined 
  \def\th@remark@my{\theorempreskipamount6\p@\@plus6\p@
    \theorempostskipamount\theorempreskipamount
    \def\theorem@headerfont{\it}\normalshape}
    \let\th@remark\th@remark@my
  \else             
    \let\o@@remark\th@remark

    \ifx\theorempostskipamount\undefined        
    \else                       
      \def\th@remark{\o@@remark
    \ifdim\theorempostskipamount < 2pt\relax
      \theorempostskipamount\theorempreskipamount
         \multiply\theorempostskipamount\tw@
         \divide\theorempostskipamount\thr@@
    \fi
      }
    \fi
  \fi
\let\myLabel\@gobble
\def\labelsONmargin{\@mparswitchfalse\def\myLabel##1{\@bsphack\marginpar
                                  {\normalshape\tiny\rm Label ##1}\@esphack}}
\ifx \url\undefined
  \def\url#1{{\tt #1}}%
\fi
\def\cyracc{\def\u##1{
                \if \i##1\char"1A%
                \else \if I##1\char"12%
                \else \accent"24 ##1\fi\fi }%
\def\"##1{\if e##1{\char"1B}%
                \else \if E##1{\char"13}%
                \else \accent"7F ##1\fi\fi }%
\def\9##1{\if##1z\char"19
\else\if##1Z\char"11
\else\if##1E\char"03
\else\if##1e\char"0B
\else\if##1u\char"18
\else\if##1U\char"10
\else\if##1A\char"17
\else\if##1a\char"1F
\else\if##1p\char"7E
\else\if##1P\char"5E
\else\if##1Q\char"5F
\else\if##1q\char"7F
\else\if##1i\char"1A
\else\if##1I\char"12
\else\if##1N\char"7D
\fi
\fi
\fi
\fi
\fi
\fi
\fi
\fi
\fi
\fi
\fi
\fi
\fi
\fi
\fi
}%
\def\cydot{{\kern0pt}}}%
\def\cydot{$\cdot$}

\ifx\Russian\undefined
        \def\Russian{0\relax
    \message{Don't know the hyphenation rules for Russian^^J
                        Please do INITeX with `input  russhyph' in the
                        command line}%
                \gdef\Russian{0\relax}%
        }
\fi
\ifx\RequirePackage\undefined           
  \def\@putname#1#2#3#4{\def\@@ref{#3}\let\old@bf\bf
        \def\bf##1{\old@bf\if?\noexpand##1?{#4}\else##1\fi}%
    #1{#2}%
        \let\bf\old@bf}
\else                       
  \def\@putname#1#2#3#4{\def\@@ref{#3}\let\old@bf\bf    
    \let\old@reset@font\reset@font          
        \def\bf##1{\old@bf\if?\noexpand##1?{#4}\else##1\fi}%
    \def\reset@font##1##2{\old@reset@font##1\if?\noexpand##2?{#4}\else##2\fi}#1{#2}%
        \let\bf\old@bf\let\reset@font\old@reset@font}
\fi
\let\my@ref=\ref
\def\ref#1{\@putname\my@ref{#1}{#1}{\tiny\rm\@@ref}}
\let\my@pageref=\pageref
\def\pageref#1{\@putname\my@pageref{#1}{#1}{\tiny\rm\@@ref}}
\let\my@cite=\cite
\def\cite#1{\@putname\my@cite{#1}{\@citeb}{\tiny\rm\@@ref}}
\makeatother
\fi
\relax
\theoremstyle{plain} 
\theorembodyfont{\sl}

\let\cal\mathcal
\let\goth\mathfrak

\def\gg{\goth g}
\def\gs{\goth s}
\def\gh{\goth h}

\def\gl{\goth l}

\def\gb{\goth b}

\def\gp{\goth p}

\def\Z{\mathbb Z}
\def\C{\mathbb C}

\DeclareMathOperator{\supp}{supp} \DeclareMathOperator{\ad}{ad}
\DeclareMathOperator{\Span}{Span} 
\DeclareMathOperator{\Hom}{Hom} 
 \DeclareMathOperator{\Id}{Id}
\DeclareMathOperator{\End}{End}\DeclareMathOperator{\im}{im}

\numberwithin{equation}{section}
\theorembodyfont{\rm}
\theoremstyle{definition}

\newtheorem{definition}{Definition}[section]

\newtheorem{example}[definition]{Example}

\theoremstyle{remark}

\newtheorem{remark}[definition]{Remark} 

\theoremstyle{plain} 
\theorembodyfont{\sl}

\newtheorem{theorem}[definition]{Theorem}
\newtheorem{lemma}[definition]{Lemma}
\newtheorem{corollary}[definition]{Corollary}
\newtheorem{proposition}[definition]{Proposition}

\begin{document}
\bibliographystyle{amsplain}
\relax

\title{ Category of $ {\mathfrak s}{\mathfrak p}\left(2n\right) $-modules with bounded weight multiplicities }

\author{Dimitar Grantcharov and Vera Serganova }

\date{ \today }

\address{ Dept. of Mathematics, San Jose State University,
San Jose, CA 95192 } \email{grantcharov\atSign{}math.sjsu.edu}

\address{ Dept. of Mathematics, University of California at Berkeley,
Berkeley, CA 94720 } \email{serganov\atSign{}math.berkeley.edu}

\maketitle

\begin{abstract}
Let ${\mathfrak g}$ be a finite dimensional simple Lie algebra.
Denote by ${\mathcal B}$ the category of all bounded weight
$\mathfrak g$-modules, i.e. those which are direct sum of their
weight spaces and have uniformly bounded weight multiplicities. A
result of Fernando shows that infinite-dimensional bounded weight
modules exist only for $\mathfrak g = \mathfrak{ sl} (n)$ and
$\mathfrak g = \mathfrak{ sp} (2n)$. If $\mathfrak g= \mathfrak{
sp} (2n)$ we show that ${\mathcal B}$ has enough projectives if
and only if $n>1$. In addition, the indecomposable projective
modules can be parameterized and described explicitly. All
indecomposable objects are described in terms of indecomposable
representations of a certain quiver with relations. This quiver is
wild for $n>2$. For $n=2$ we describe all indecomposables  by
relating the blocks of $\mathcal B$ to the representations of the
affine quiver $A_3^{(1)}$.
\end{abstract}

\section{Introduction}

To classify all indecomposable objects in a category of
representations is usually a challenging and difficult problem. It
is often the case that there are not enough projectives or the
category itself is wild. A classical example of a wild category
with enough projectives is the category ${\mathcal O}$ introduced
by Bernstein-Gelfand-Gelfand in 1967. The simple objects in this
category are highest weight modules, and the indecomposable
projectives are described by the celebrated BGG reciprocity law.

A natural generalization of the category $\mathcal O$ is the
category of all weight (not necessarily highest weight) modules.
Weight modules have attracted considerable mathematical attention
in the last 20 years and appeared in works of G. Benkart, D.
Britten, S. Fernando, V. Futorny, and F. Lemire, \cite{BBL},
\cite{BL1}, \cite{BL2}, \cite{F}, \cite{Fu}.
 A major breakthrough was the recent
classification of O. Mathieu, \cite{M}, of all simple weight
modules with finite weight multiplicities over finite dimensional
reductive Lie algebras. A crucial role in this classification is
played by the category ${\mathcal B}$ of bounded weight modules,
i.e. those for which the set of weight multiplicities is uniformly
bounded. This is due to the fact that, as Fernando showed in
\cite{F}, every simple weight module $M$ with finite weight
multiplicities is obtained by a parabolic induction from a simple
module $S$ in ${\mathcal B}$ (in fact $S$ has equal weight
multiplicities). An important observation of Mathieu is that the
direct sum of all simple objects in a single block of $\mathcal B$
form a so-called coherent family which is parameterized by a
highest weight module, i.e. an object in ${\mathcal O}$.

In the present paper we initiate a study of the category of
bounded modules. A result of Fernando shows that infinite
dimensional bounded weight modules exist only for Lie algebras of
type $A$ and $C$ (\cite{F}, \cite{M}). As a first step in our
project we consider the Lie algebra $\gg = \gs \gp (2n)$. This
case is simpler in terms of the classification of Mathieu as a
semisimple irreducible coherent family over $\gs \gp (2n)$ is
determined uniquely by its central character. The case of $\gs \gl
(n+1)$ is more delicate and one has to consider three separate
cases for the central character: regular integral, singular, and
nonintegral.

One of the main results in the paper is providing a complete
classification of all indecomposable projective objects in
${\mathcal B}$. An interesting observation is that if $n=1$, i.e.
$\gg =\gs \gl (2)$, the category ${\mathcal B} $ does not contain
any projective objects. The picture is totally different for the
higher dimensional algebras as for $n>1$ each simple object has a
projective cover.

In order to describe the indecomposable objects of $\mathcal B$ we
first show that this category is equivalent to the category of
weight modules over the Weyl algebra ${\mathcal A}_n$ (see Lemma
\ref{lm7} and Corollary \ref{cor3}). We then conclude that each
block ${\mathcal B}^{\chi}$ of $\mathcal B$ is equivalent to the
category of a certain quiver with relations. This quiver is wild
if and only if $n>2$. In the case $n=2$ indecomposable
representations of the quiver can be expressed in terms of the
affine quiver $A_3^{(1)}$, the theory of which is well
established. In addition, in section \ref{explic} we provide an
explicit description of all indecomposable bounded modules over
$\gs \gp (4)$ in terms of the twisted localization correspondence.

We show also that there are not enough projectives in the category
of all weight $\gg$-modules with finite weight multiplicities (see
Example \ref{exproj}) which provides an additional motivation to
focus our attention on the bounded modules only.

\section{Weight modules over the Weyl algebra}

The ground field is $\mathbb C$. By ${\mathcal A}_{n}$ we denote
the Weyl algebra, i.e. the algebra of polynomial differential
operators on $ {\mathbb A}^{n} $. Let $ t_{1},\dots
,t_{n},\partial_{1},\dots ,\partial_{n} $ be the standard
generators of $ {\mathcal A}_{n} $. Recall that the following
relations hold
\begin{equation}
\left[t_{i},\partial_{j}\right]=\delta_{ij}\text{,
}\left[t_{i},t_{j}\right]=\left[\partial_{i},\partial_{j}\right]=0.
\notag\end{equation} In what follows we will consider ${\mathcal
A}_{n}$ as a Lie algebra over $\mathbb C$. Let $M$ be an
${\mathcal A}_{n}$-module. We say that $ M $ is a {\it weight
module} if
\begin{equation}
M=\bigoplus_{\mu\in{\mathbb C}^{n}} M^{\mu}, \notag
\end{equation}
where $ M^{\mu}:=\left\{m\in M \mid
t_{i}\partial_{i}\left(m\right)=\mu_{i} m \mbox{ for all }
i\right\} $ and $\mu = (\mu_1,...,\mu_n) \in \C^n$. The space  $
M^{\mu}$ is the {\it weight space of weight $\mu$} and $\dim
M^{\mu}$ is the {\it weight multiplicity} of $ M^{\mu}$. We say
that  $ M $ is {\it multiplicity free} if $ \dim M^{\mu} \leq 1$.
The {\it support} of $M$ is the set $\operatorname{supp}
M:=\left\{\mu\in{\mathbb C}^{n} \mid M^{\mu}\not=0\right\} $.

The Lie algebra $ {\mathcal A}_{n}$ acts on itself via the adjoint
map $\ad :  {\mathcal A}_{n} \to \End ( {\mathcal A}_{n}),
\ad(x)(y):=[x,y]$. The elements $ t_{1}\partial_{1},\dots
,t_{n}\partial_{n} $ act diagonally on $ {\mathcal A}_{n} $. The
adjoint action induces a $\Z^n$-grading of $ {\mathcal A}_{n} $
via the root decomposition:
\begin{equation}
{\mathcal A}_{n}=\bigoplus_{\alpha\in P} {\mathcal
A}_{n}^{\alpha}, \notag\end{equation} where $P={\mathbb Z}^{n}$ is
considered as a sublattice of $ {\mathbb C}^{n} $ with the
standard generators $\varepsilon_1,\dots, \varepsilon_n$. The
following lemma follows by a direct verification.

\begin{lemma} \label{lm0}\myLabel{lm0}\relax  $ {\mathcal A}_{n}^{0} $ is a free
commutative algebra with generators
$ t_{1}\partial_{1},\dots ,t_{n}\partial_{n} $ and each $ {\mathcal A}_{n}^{\alpha} $ is a free left $ {\mathcal A}_{n}^{0} $-module of rank 1.

\end{lemma}

\begin{example} \label{ex1}\myLabel{ex1}\relax
Let $ \mu\in{\mathbb C}^{n} $, and let $ t^{\mu} $ stand for $
t_{1}^{\mu_{1}}\dots t_{n}^{\mu_{n}} $. The vector space $
F\left(\mu\right)=t^{\mu}{\mathbb C}\left[t_{1}^{\pm1},\dots
,t_{n}^{\pm1}\right] $ has a natural structure of an $ {\mathcal
A}_{n} $-module. It is an easy exercise to check that $
F\left(\mu\right) $ is a multiplicity free $ {\mathcal A}_{n}
$-module with $ \operatorname{supp} F\left(\mu\right)=\mu+P $.

\end{example}

\begin{lemma} \label{lm1}\myLabel{lm1} \relax  The $ {\mathcal A}_{n} $-module
$ F\left(\mu\right) $ is indecomposable. It is irreducible if and
only if $ \mu_{i}\notin {\mathbb Z} $ for all $ i=1,\dots ,n $.

\end{lemma}

\begin{proof} Suppose that $F(\mu)=M_1 \oplus M_2$. Since $F(\mu)$ is multiplicity free,
$\operatorname {supp} F (\mu)$ is a disjoint union of
$\operatorname {supp} M_1$ and $\operatorname {supp} M_2$.
Therefore one can find $p \in \{1, 2\}$, $\nu \in \operatorname
{supp} M_p$, and $i \in \{ 1, 2, ..., n\} $ such that $\nu
+\varepsilon _i \notin \operatorname {supp} M_p$. Then $t_i v=0$
whenever $v\in F(\mu)^{\nu}$. This is impossible because $F(\mu)$
is free over $\mathbb C \left[ t_1,\dots, t_n\right]$.

To prove the second statement we first assume that $\mu _i \in
\mathbb Z$ for some $i$. Since $F(\mu)$ is isomorphic to $F(\mu +
\gamma)$ for any $\gamma \in P$ we may assume that $\mu _i =0$.
Then one easily checks that $t^{\mu}\mathbb C
\left[t_{1}^{\pm1},\dots,t_{i-1}^{\pm1},t_i, t_{i+1}^{\pm1},
\dots, t_{n}^{\pm1}\right]$ is a submodule of $F(\mu)$. Finally,
if $\mu_i \notin \mathbb Z$ for all $i$, any element of
$F(\mu)^{\nu}$ generates $F(\mu)$. Hence $F(\mu)$ is irreducible.
\end{proof}

Denote by $ {\mathcal F}_{n} $ the category of weight $ {\mathcal
A}_{n} $-modules with finite weight multiplicities.

\begin{lemma} \label{lm2}\myLabel{lm2}\relax  The category
$ {\mathcal F}_{n} $ splits into a direct sum of blocks
\begin{equation}
\bigoplus_{\bar{\nu}\in{\mathbb C}^{n}/P}{\mathcal
F}_{n}^{\bar{\nu}}, \notag\end{equation} where the sum runs over
all distinct classes $\bar{\nu}:= \nu + P$ in $\C^n/P$ and $
{\mathcal F}_{n}^{\bar{\nu}} $ is the subcategory of all modules $
M $ such that $ \operatorname{supp} M\subset \bar{\nu}$.
\end{lemma}

\begin{proof} Let $M \in {\mathcal F}_n$. For any $\bar{\nu} \in {\mathbb C}^{n}/P$ let
\begin{equation}
M(\bar{\nu}):=\bigoplus_{\mu \in \bar{\nu}} M^{\mu}.
\notag\end{equation} Obviously, $M(\bar{\nu})$ is a submodule of
$M$ and
\begin{equation}
M=\bigoplus_{\bar{\nu}\in{\mathbb C}^n/P} M(\bar{\nu}).
\notag\end{equation} This proves the lemma.
\end{proof}

For each $ \mu\in P $ put
\begin{equation}
P\left(\mu\right):={\mathcal A}_{n}\otimes_{{\mathcal A}_{n}^{0}}
C_{\mu}, \notag\end{equation} where $ C_{\mu} $ denotes the unique
$ 1 $-dimensional $ {\mathcal A}_{n}^{0} $-module of weight $ \mu
$.

\begin{theorem} \label{th1} \myLabel{th1}\relax
\begin{enumerate}
\item $ P\left(\mu\right) $ is a multiplicity free module with $
\supp P\left(\mu\right)=\mu+P $;

\item $ P\left(\mu\right) $ has a unique irreducible quotient
which we denote by $ L\left(\mu\right) $;

\item If $ M $ is an irreducible module in $ {\mathcal F}_n $ such
that $ \mu\in\operatorname{supp} M $, then $ M $ is isomorphic to
$ L\left(\mu\right) $;

\item $ P\left(\mu\right) $ is indecomposable;

\item $ P\left(\mu\right) $ is a projective module in the category
$ {\mathcal F}_n $;

\item Every indecomposable projective module in the category $
{\mathcal F}_n $ is isomorphic to $ P\left(\mu\right)$ for some
$\mu$.
\end{enumerate}
\end{theorem}

\begin{proof} The first statement follows from Lemma~\ref{lm0}. To show (2) it suffices
to prove that $ P\left(\mu\right) $ has a unique maximal proper submodule. Indeed, $ N $ is a
proper submodule of $ P\left(\mu\right) $ iff $ \mu\notin\operatorname{supp} N $. Since
\begin{equation}
\operatorname{supp} \left(N_{1}\oplus
N_{2}\right)=\operatorname{supp} N_{1} \cup \operatorname{supp}
N_{2}, \notag\end{equation} the sum of all proper submodules of $
P\left(\mu\right) $ is proper. (3) follows from the Frobenius
reciprocity theorem, and (4) follows from (2). To prove (5)
consider an exact sequence
\begin{equation}
0 \to N \to S \xrightarrow[]{p} P\left(\mu\right) \to\text{ 0.}
\notag\end{equation}
Since the sequence
\begin{equation}
0 \to N^{\mu} \to S^{\mu} \to P\left(\mu\right)^{\mu} \to\text{ 0}
\notag
\end{equation}
of $ {\mathcal A}_{n}^{0} $-modules splits,  there is a map $
i\colon P\left(\mu\right)^{\mu}\cong C_{\mu} \to S^{\mu} $ such
that $ i\circ p=id. $ By the Frobenius reciprocity theorem, the
map $ i $ induces a map $ j:P\left(\mu\right) \to S $ for which $
j\circ p= \Id. $ Hence $ P\left(\mu\right) $ is projective. To
prove (6) let $S$ be an indecomposable projective module. Then we
have a surjective map $a: S \to L(\mu)$ for some irreducible
module $L(\mu)$. Let $r: P(\mu) \to L(\mu)$ be the canonical map.
Then there exist $b: S \to P(\mu)$ and $c:P(\mu) \to S$ such that
$a \circ c=r$ and $r \circ b=a $. Then $r \circ  b \circ c=r$, and
therefore $b \circ c \neq 0$. On the other hand, one can easily
see that $\operatorname{End}_{\gg}P(\mu)=\mathbb C$. Hence $b
\circ c$ is an automorphism. In particular, $b$ is surjective.
Then $P(\mu)$ is isomorphic to a direct summand of $S$. But $S$ is
indecomposable, so $S$ is isomorphic to $P(\mu)$.
\end{proof}

\begin{corollary} \label{cor111}\myLabel{cor111}\relax
Let $M$ and $N$ be simple modules in $\mathcal F_n$. Then $M$ and
$N$ are nonisomorphic  if and only if $\operatorname {supp} M$ and
$\operatorname {supp} N$ are disjoint.
\end{corollary}

\begin{proof} If $\mu \in \operatorname {supp} M \cap \operatorname {supp} N$,
then both modules are quotients of $P(\mu)$. By Theorem \ref{th1},
(2), $P(\mu)$ has a unique simple quotient, and thus $M$ and $N$
are isomorphic.
\end{proof}

Let $M\in \mathcal F_n$ and $M=\oplus _{\mu \in \operatorname
{supp} M}M^{\mu}$. Set $M^*:=\oplus _{\mu \in \operatorname {supp}
M}(M^{\mu})^*$. Define the action of $\mathcal A_n$ on $M^*$ by
\begin{equation}
\partial_i \cdot \tau (v)=\tau (t_i \cdot v), t_i\cdot \tau (v)=\tau
(\partial_i \cdot v) \notag\end{equation} for any $v\in M, \tau
\in M^*$. It is easy to check that $M^*\in \mathcal F_n$ and
$\operatorname {supp} M=\operatorname {supp} M^*$. Moreover, $*$
is an exact contravariant functor on $\mathcal F_n$ which maps
projective objects to injective ones and preserves the simple
objects.

To obtain a complete description of all irreducible and
indecomposable projectives in each block $ {\mathcal
F}_{n}^{\bar{\nu}} $ we observe that
\begin{equation}
{\mathcal A}_{n}\cong{\mathcal A}_{1}\otimes{\mathcal
A}_{1}\otimes\dots \otimes{\mathcal A}_{1}. \notag\end{equation}
Therefore every irreducible object in $ {\mathcal F}_{n} $ is a
tensor product of irreducibles in $ {\mathcal F}_{1} $, and by
Theorem~\ref{th1}, the same holds for the indecomposable
projectives. Hence, it is enough to describe the blocks of $
{\mathcal F}_{1} $. This description is obtained in the following
lemma.

\begin{lemma} \label{lm3}\myLabel{lm3}\relax For any $ \bar{\nu} \not= \bar{0} $,
the block $ {\mathcal F}_{1}^{\bar{\nu}} $ is semi-simple and has
exactly one up to isomorphism irreducible object $
F\left(\mu\right) $, $ \mu\in \bar{\nu} $. The block $ {\mathcal
F}_{1}^{0} $ has two isomorphism classes of simple objects: $
L\left(0\right) $ and $ L\left(-1\right) $. The structure of the
indecomposable projective modules is described by the following
exact sequences
\begin{equation}
0 \to L\left(-1\right) \to P\left(0\right) \to L\left(0\right) \to
0, \; 0 \to L\left(0\right) \to P\left(-1\right) \to
L\left(-1\right) \to 0. \notag\end{equation}
\end{lemma}

\begin{proof} By Lemma \ref{lm1} $F(\mu)$ is irreducible iff
$\mu=\mu _1\notin \mathbb Z$. Clearly, in this case $F(\mu)$ is
isomorphic to $P(\mu)$, therefore $\mathcal F_1^{\bar{\mu}}$
contains one up to an isomorphism indecomposable object $F(\mu)$
which is both projective and simple.

If $\bar{\nu} =\bar{0}$, then $F(0) \cong F(n)$ for any $n \in
\mathbb Z$ and a simple calculation leads to the exact sequence
\begin{equation}
0 \to L\left(0\right) \to F\left(0\right) \to L\left(-1\right) \to
0. \notag\end{equation} The Frobenius reciprocity implies that
there is a surjective homomorphism $P(-1) \to F(0)$, which is an
isomorphism because both modules are multiplicity free and have
the same support. By Corollary \ref{cor111} every simple object in
${\mathcal F}_{1}^{\bar{0}}$ is a subquotient of $F(0)$. Finally,
by similar arguments $P(0)\cong F(0)^*$, which leads to the exact
sequence for $P(0)$.
\end{proof}

\begin{remark} \label{rem1}\myLabel{rem1}\relax  One can use also
the following geometric description. $ L\left(0\right) $ is
isomorphic to $ {\mathbb C}\left[t\right] $, $ P\left(-1\right) $
is isomorphic to $ {\mathbb C}\left[t,t^{-1}\right] $, and $
L\left(-1\right) $ is a module generated by the $ \delta
$-function concentrated at zero on $ {\mathbb C}^{1} $.

\end{remark}

\begin{corollary} \label{cor1}\myLabel{cor1}\relax
Let $ \nu\in{\mathbb C}^{n}$ and $I(\bar{\nu}):=\left\{i\leq n
\mid \nu_{i} \in \Z\right\} $. Then all indecomposable projective
modules and all irreducible modules of $ {\mathcal
F}_{n}^{\bar{\nu}} $ are parameterized by the set $\mathcal S$ of
all maps $ s : I(\bar{\nu}) \to \left\{0,-1\right\} $. More
precisely, $ P\left(s\right) $ is the tensor product of $
P\left(\nu_{j}\right) $ for $ j\notin I(\bar{\nu}) $ and $
P\left(s\left(i\right)\right) $ for $ i\in I(\bar{\nu}) $. The
same description works for the irreducibles.

\end{corollary}

Since $ {\mathcal F}_{n}^{\bar{\nu}} $ has finitely many
irreducible modules and each irreducible has a unique
indecomposable projective cover, the category $ {\mathcal
F}_{n}^{\bar{\nu}} $ is equivalent to the category of
finite-dimensional $ E^{\nu} $-modules, where
\begin{equation}
E^{\nu}:=\operatorname{End}_{{\mathcal A}_{n}}\left(\oplus_{s \in
{\mathcal S}}P\left(s\right)\right). \notag\end{equation}
Furthermore,
\begin{equation}
E^{\nu}\cong E^{\nu_{1}}\otimes\dots \otimes E^{\nu_{n}}.
\label{equ0}\end{equation}\myLabel{equ0,}\relax Observe that $
E^{\nu_{i}}\cong{\mathbb C} $ whenever $ \nu_{i}\notin \Z$. Let $
V_{1} $ be the quiver
\begin{diagram}
\bullet & \pile{\rTo^{\varphi^+} \\ \lTo_{\varphi^{-}}} & \bullet
 \end{diagram}
with relations $ \varphi^{+}\varphi^{-}=\varphi^{-}\varphi^{+}=0
$. Then one can see easily that $ E^{\nu_{i}}\cong{\mathbb
C}\left(V_{1}\right) $ in the case $ \nu_{i} \in \Z$.

Define the quiver $ V_{k} $ in the following way. The vertices of
$ V_{k} $ are the vertices of the cube in $ {\mathbb R}^{k} $ with
coordinates $1$ or $-1$. The edges are the edges of the cube with
two possible orientations. We call a path on the cube {\it
admissible} if each coordinate function is weakly monotonic along
the path. Finally, we impose the following relations: each
non-admissible path is zero, every two admissible paths with the
same start and the same end points are equal.

\begin{theorem} \label{th2}\myLabel{th2}\relax
Let $\bar{\nu}\in{\mathbb C}^{n}/P $ and $ k $ be the number of
all $ i $ for which $ \bar{\nu}_{i}=\bar{0} $. Then $ {\mathcal
F}_{n}^{\bar{\nu}}$ is equivalent to the category of
representation of the quiver $ V_{k} $.

\end{theorem}

\begin{lemma} \label{lm4}\myLabel{lm4}\relax  For $ k\geq3 $ the quiver $ V_{k} $ is wild.

\end{lemma}
\begin{proof} Choose a subquiver $V_3 \subset V_k$ in an arbitrary way. Then
  choose $W_3 \subset V_3$ to be a maximal subquiver without
cycles. Every representation of $W_3$ can be extended to a
representation of  $V_k$ trivially: every arrow of $V_k$ which is
not in $W_3$ is represented by the zero map. Since $W_3$ is wild,
$V_k$ is wild as well.
\end{proof}

The indecomposable representations of $ V_{1} $ are easy to describe.

\begin{lemma} \label{lm5}\myLabel{lm5}\relax
The quiver $ V_{1} $ has four isomorphism classes of
indecomposable representations with dimension functions $(1,0)$,
$(0,1)$, $(1,1)$, and $(1,1)$, respectively.

\end{lemma}

\begin{proof} Consider an indecomposable representation of $ V_{1} $. Let $ A_{1} $ and
$ A_{2} $ be the spaces attached to the vertices of $ V_{1} $, and
let
\begin{equation}
\varphi^{+}:A_{1} \to A_{2}\text{, }\varphi^{-}:A_{2} \to A_{1}
\notag\end{equation}
 be the corresponding maps. We have that
$\varphi^{+}\varphi^{-}=\varphi^{-}\varphi^{+}=0 $. Choose $
B_{1}\subset A_{1} $ and $ B_{2}\subset A_{2} $, so that $
A_{1}=B_{1}\oplus\operatorname{Ker} \varphi^{+} $ and $
A_{2}=B_{2}\oplus\operatorname{Ker} \varphi^{-} $. Then the
representation splits into the direct sum
\begin{equation}
(\varphi^{+}:B_{1} \to \operatorname{Ker} \varphi^{-}) \oplus
(\varphi^{-}:B_{2} \to \operatorname{Ker} \varphi^{+}).
\notag\end{equation} Thus either $ B_{1}=0 $ or $ B_{2}=0 $, and
the problem is reduced to the quiver
\begin{diagram}
\bullet & \rTo & \bullet \end{diagram} which is
well-understood.\end{proof}

To describe the indecomposable representations of $ V_{2} $ we
first introduce some notation. By $ \rho_{1} $ we denote the
following indecomposable representation of $ V_{2} $
\begin{equation}
\begin{CD}
 {\mathbb C} @>>> {\mathbb C}
\\
@VVV @VVV
\\
{\mathbb C} @>>> {\mathbb C}
\\
\end{CD}
\notag\end{equation} where all the arrows are represented by the
identity maps and all inverse arrows are represented by the zero
maps. One obtains $ \rho_{2} $, $ \rho_{3} $, $ \rho_{4} $ from $
\rho_{1} $ by rotating the picture by $ 90^{\circ} $ one, two, or
three times, respectively.

We next introduce the quivers $ A $ and $ B $:
\begin{equation}
\begin{CD}
 A_{11} @>>> A_{12}
\\
@VVV @AAA
\\
A_{21} @<<< A_{22}
\\
\end{CD}
\notag\end{equation}
\begin{equation}
\begin{CD}
 B_{11} @<<< B_{12}
\\
@AAA @VVV
\\
B_{21} @>>> B_{22}
\\
\end{CD}
\notag\end{equation}

\medskip
Any indecomposable representation of $ A $ or $ B $ induces an
indecomposable representation of $ V_{2} $ if we represent all
reverse arrows in $ V_{2} $ by the zero maps.

\begin{lemma} \label{lm6}\myLabel{lm6}\relax  Any indecomposable representation of $ V_{2} $ is either
isomorphic to one of the representations $
\rho_{1},\rho_{2},\rho_{3},\rho_{4} $ or induced by an
indecomposable representation of $ A $ or $ B $.

\end{lemma}

\begin{proof}
Consider some indecomposable indecomposable representation $ \rho
$ of $ V_{2} $:

\begin{diagram}
C_{11} & \pile{\rTo^{\varphi^{+}} \\ \lTo_{\varphi{-}}} & C_{12} \\
\dTo^{\xi^{+}} \uTo_{\xi^{-}} & & \uTo^{\eta^{-}} \dTo_{\eta^{+}} \\
C_{21} & \pile{\lTo^{\psi^{-}} \\ \rTo_{\varphi{+}}} & C_{22}
\end{diagram}
Assume that there is $ v\in C_{11} $ such that $
\eta^{+}\varphi^{+}\left(v\right)\not=0 $. The relations of $V_2$
imply that
\begin{equation}
\psi^{+}\xi^{+}\left(v\right)=\eta^{+}\varphi^{+}\left(v\right).
\notag\end{equation} One can see easily that $ v $ generates a
subrepresentation $\rho'$ of $ \rho $ isomorphic to $ \rho_{1} $.
Moreover, $\rho'$ is a direct summand of $\rho$, since each of the
vectors $v$, $ \varphi^{+}\left(v\right) $, $
\eta^{+}\varphi^{+}\left(v\right) $ and $ \xi^{+}\left(v\right) $
does not belong to the sum of images of all reverse maps, i.e.
$$
v \notin \im(\xi^{-}) + \im(\phi^{+}), \, \varphi^+ (v) \notin \im
(\eta^{-1}), \, \xi^{+}(v)\notin \im (\varphi^{-}).
$$
Indeed, say $ v=\xi^{-}\left(u\right)+\varphi^{-}\left(w\right) $.
Then
\begin{equation}
\eta^{+}\varphi^{+}\left(v\right)=\eta^{+}\varphi^{+}\xi^{-}\left(u\right)\not=0,
\notag
\end{equation}
which contradicts the relations. Thus in this case $
\rho\cong\rho_{1} $. In the same way, if we start with $ v\in
C_{12} $ we will conclude that $ \rho\cong\rho_{2}$, etc.

Let us assume now that $ \rho $ is not isomorphic to $
\rho_{1},\rho_{2},\rho_{3} $ or $ \rho_{4} $. Then the above
argument shows that a composition of any two arrows is the zero
map. Let $ U_{ij} $ be the intersection of the kernels of the two
maps starting at $ C_{ij} $. Write $ C_{ij}=U_{ij}\oplus D_{ij} $
choosing $ D_{ij} $ in an arbitrary way. Then $
\rho=\pi_{1}\oplus\pi_{2} $, where $ \pi_{1} $ is the following
representation of $ A $
\begin{equation}
\begin{CD}
 D_{11} @>>> U_{12}
\\
@VVV @AAA
\\
U_{21} @<<< D_{22}
\\
\end{CD}
\notag\end{equation}
and $ \pi_{2} $ is the following representation of $ B $
\begin{equation}
\begin{CD}
 U_{11} @<<< D_{12}
\\
@AAA @VVV
\\
D_{21} @>>> U_{22}
\\
\end{CD}
\notag\end{equation}

Hence the lemma is proved.\end{proof}

Since the quivers $ A $ and $ B $ are isomorphic to affine Dynkin graph $ A_{3}^{\left(1\right)} $,
we can use the general theory of representation of tame quivers.

\section{The algebra $ {\protect \mathcal A}_{n}^{ev} $ }

Let $ Q $ be a sublattice of index $2$ in $ P={\mathbb Z}^{n} $
consisting of all $ \left(\mu_{1},\dots ,\mu_{n}\right) $ such
that $ \mu_{1}+\dots +\mu_{n}\in2{\mathbb Z} $. Define
\begin{equation}
{\mathcal A}_{n}^{ev}=\bigoplus_{\alpha\in Q}{\mathcal
A}_{n}^{\alpha}. \notag\end{equation} Clearly, $ {\mathcal
A}_{n}^{ev} $ is a Lie subalgebra of $ {\mathcal A}_{n} $.

Denote by $ {\mathcal F}_{n}^{ev} $ the category of weight $
{\mathcal A}_{n}^{ev} $-modules with finite weight multiplicities.
As in Lemma~\ref{lm2} one has a block decomposition
\begin{equation}
{\mathcal F}_{n}^{ev}=\bigoplus_{\bar{\nu}\in{\mathbb
C}^{n}/Q}\left({\mathcal F}_{n}^{ev}\right)^{\bar{\nu}}.
\notag\end{equation}

Let $ \bar{\nu} \in{\mathbb C}^{n}/Q $, and let $ \bar{\mu}
\in{\mathbb C}^{n}/P $ be the image of $ \bar{\nu} $ under the the
natural projection $ {\mathbb C}^{n}/Q \to {\mathbb C}^{n}/P$.
Define two functors
\begin{equation}
\operatorname{Ind}:\left({\mathcal F}_{n}^{ev}\right)^{\bar{\nu}}
\to {\mathcal F}_{n}^{\bar{\mu}}\text{, }\operatorname{Res}\colon
{\mathcal F}_{n}^{\bar{\mu}} \to \left({\mathcal
F}_{n}^{ev}\right)^{\bar{\nu}} \notag
\end{equation}
by putting
\begin{equation}
\operatorname{Ind} \left(M\right)={\mathcal
A}_{n}\otimes_{{\mathcal A}_{n}^{ev}}M\text{, }\operatorname{Res}
\left(N\right)=\oplus_{\gamma\in \bar{\nu}} N^{\gamma}.
\notag
\end{equation}
The following lemma is straightforward.

\begin{lemma} \label{lm7}\myLabel{lm7}\relax
The functors $ \operatorname{Ind} $ and $ \operatorname{Res} $
establish an equivalence of the categories $ \left({\mathcal
F}_{n}^{ev}\right)^{\bar{\nu}} $ and $ {\mathcal
F}_{n}^{\bar{\mu}} $.
\end{lemma}

\section{Twisted localization of bounded modules} \label{twloc}

Let $\gg = \gs \gp (2n)$ or $\gg = \gs \gl (n+1)$, and $\gh$ be a
Cartan subalgebra of $\gg$. Let  $\Delta = \Delta (\gg, \gh)$ be
the root system, and $Q$ be the root lattice of $\gg$.  For every
$\alpha \in \Delta$ fix a standard triple $\{ e_{\alpha},
f_{\alpha}, h_{\alpha}\}$ such that  $e_{\alpha} \in \gg^{\alpha},
f_{\alpha} \in \gg^{-\alpha}$ and $[e_{\alpha}, f_{\alpha}] =
h_{\alpha}$. Let $U:= U(\gg)$ be the universal enveloping algebra
of $\gg$, $Z:= Z(\gg)$ be its center, and $Z':=\Hom(Z, \C)$.

By $ {\mathcal B}^{\chi} $ we denote the category of weight $
{\mathfrak g} $-modules with bounded weight multiplicities
admitting generalized central character $ \chi \in Z'$. In other
words, $ M\in{\mathcal B}^{\chi} $ if
\begin{equation}
M=\bigoplus_{\mu\in{\mathfrak h}^{*}} M^{\mu},
\notag
\end{equation}
there exists $ C_{M}$ such that $\dim  M^{\mu}<C_{M} $ for all
$\mu \in {\mathfrak h}^*$, and for each $ m\in M $ and $z \in Z$
there exists $ N $ such that
\begin{equation}
\left(z-\chi\left(z\right)\right)^{N}m=0. \notag\end{equation}

Put ${\cal B}:= \cup_{\chi \in Z'} {\cal B}^{\chi}$. In what
follows we assume that all $\gg$-modules are {\it
bounded}\footnote{In \cite{M} Mathieu uses the term ``admissible''
weight module, but to avoid confusion with Harish-Chandra modules
of finite type we prefer to use the term ``bounded'' weight
module}, i.e. in ${\cal B}$. Following the approach in \cite{M},
we  recall some facts about the localization of (bounded) weight
modules with respect to a
 set of commuting roots. Let $\Gamma =\{ \gamma_1,..., \gamma_l\}
\subset \Delta$ be a linearly independent subset of $Q$ for which
$\gamma_i + \gamma_j \notin \Delta$. The set $\{ f_{\gamma_1},...,
f_{\gamma_l}\}$ generates a multiplicative subset  $F_{\Gamma}$ of
$U$ which satisfies Ore's localizability conditions. Let
$U_{F_\Gamma}$ be the localization of $U$ relative to
$F_{\Gamma}$.

A $\gg$-module $M$ is called $\Gamma$-{\it injective} (
$\Gamma$-{\it bijective}) if $f_{\gamma}$ acts injectively
(bijectively) on $M$ for every $\gamma$ in $\Gamma$. For any
$\gg$-module $M$ we define the $\Gamma${\it -localization} ${\cal
D}_{\Gamma}M$ of $M$ by  ${\cal D}_{\Gamma}M:= U_{F_\Gamma}
\otimes_U M$. If $M$ is $\Gamma$-injective, then $M \subset {\cal
D}_{\Gamma}M$.  Note that if $\Gamma = \Gamma_1 \cup \Gamma_2$ we
have ${\cal D}_{\Gamma_1} {\cal D}_{\Gamma_2} = {\cal
D}_{\Gamma_2} {\cal D}_{\Gamma_1} = {\cal D}_{\Gamma}$ over the
set of all $\Gamma$-injective modules.

\begin{example} \label{ex_loc}
Let $\gg = \gs \gp(2n)$, $\gb$ be the standard Borel subalgebra
with basis $\{ \varepsilon_1 - \varepsilon_2,...,\varepsilon_{n-1}
- \varepsilon_n, 2\varepsilon_n\}$, and $\Gamma :=
\{2\varepsilon_1,...,2\varepsilon_n \}$. Then every simple
$\gb$-highest weight module $M = L_B(\lambda)$ is
$\Gamma$-injective. Furthermore, if $M$ is bounded, then ${\cal
D}_{\Gamma}M$ has $2^n$ simple subquotients all of which are
highest weight modules (with respect to different Borel
subalgebras). This is proved in \cite{BKLM} and a detailed
description of ${\cal D}_{\Gamma}M$ for $\gg =\gs \gp(4)$ will be
provided in section \ref{explic}.
\end{example}

The preceding example is a part of more general picture which is
summarized in the following statement. The proof uses a
combinations of statements (Lemma 4.5, Proposition 4.8, and Lemma
9.2) in \cite{M} and is based on
 Mathieu's description of the coherent extensions
of  bounded $\gs \gp (2n)$-modules.

\begin{proposition} \label{simple}
Let $\gg = \gs \gp (2n)$ and $M$ be a simple module in ${\cal
B}^{\chi}$. There is a subset $\Gamma$ of $\Delta$ consisting of
$n$ long roots for which $M$ is $\Gamma$-injective. The set of all
simple subquotients of ${\cal D}_{\Gamma} M$ coincides with the
set of all simple modules $N$ in  ${\cal B}^{\chi}$ for which
$\supp N \subset \supp {\cal D}_{\Gamma} M = \supp M + Q$.
\end{proposition}

Recall now the definition of a generalized conjugation in
$U_{F_{\Gamma}}$ introduced in \cite{M}. Let $\mu =x_1\gamma_1 +
... + x_l \gamma_l \in \Span_{\C}\Gamma \subseteq \gh^*$. For $u
\in U_{F_{\Gamma}}$, $v \in N$ set
$$
\Theta_{(x_{1},\dots,x_{l})}(u):= \sum\limits_{0\leq
i_{1},\dots,i_{l}\leq N(u)} (^{x_{1}}_{i_{1}})\dots
(^{x_{l}}_{i_{l}})\,\mbox{ad}(f_{\gamma_1})^{i_{1}}\dots
\mbox{ad}(f_{\gamma_l})^{i_{l}}(u) \,f_{\gamma_1}^{-i_{1}}\dots
f_{\gamma_l}^{-i_{l}},
$$
 where $(^{x}_{i}) :=
x(x-1)...(x-i+1)/i!$ for $x \in \C$ and $i \in \Z_+ \cup \{ 0\} $.
Note that for $(x_1,...,x_l) \in \Z^l$ we have
$\Theta_{(x_{1},\dots,x_{l})}(u) =
f_{\gamma_1}^{x_1}...f_{\gamma_l}^{x_l} u
f_{\gamma_1}^{-x_1}...f_{\gamma_l}^{-x_l}$. For a
$U_{F_{\Gamma}}$-module $N$ by $\Phi^{\mu}_{\Gamma} N$ we denote
the $U_{F_{\Gamma}}$-module $N$ twisted by the action
 $$
u \cdot v^{\mu} :=
 ( \Theta_{(x_{1},\dots,x_{l})}(u)\cdot v)^{\mu},
$$
where $u \in U_{F_{\Gamma}}$, $v \in N$, and $v^{\mu}$
stands for the element $v$ considered as an element of
$\Phi^{\mu}_{\Gamma} N$. In particular, $v^{\mu} \in N^{\lambda +
\mu}$ whenever $v \in N^{\lambda}$. The following lemma is
straightforward.

\begin{lemma} \label{lmnew}
$(i)$ $\Phi^{\mu}_{\Gamma} \circ \Phi^{\nu}_{\Gamma} =\Phi^{\mu
+\nu}_{\Gamma} $, in particular, $\Phi^{\mu}_{\Gamma} \circ
\Phi^{-\mu}_{\Gamma} =\operatorname {Id}$;

$(ii)$ $\Phi^{\mu}_{\Gamma} = \Id$ whenever $\mu \in Q$;

$(iii)$ $M$ is an indecomposable $U_{F_{\Gamma}}$-module if and
only
 if $\Phi^{\mu}_{\Gamma} M$ is indecomposable.
\end{lemma}

For a $\Gamma$-injective module $M$ and $\mu \in \gh^*$ we define
the {\it twisted localization ${\cal D}_{\Gamma}^{\mu} M$ of $M$
relative to $\Gamma$ and $\mu$} by ${\cal D}_{\Gamma}^{\mu} M:=
\Phi^{\mu}_{\Gamma} {\cal D}_{\Gamma} M $. The twisted
localization plays a major role in the theory of coherent families
introduced by Mathieu. An example of such family is the {\it
coherent extension} ${\cal E}(M):=\oplus_{\bar{\mu} \in \gh^*/Q}
{\cal D}_{\Gamma}^{\bar{\mu}} M$ of $M$. Here ${\cal
D}_{\Gamma}^{\bar{\mu}} M:= {\cal D}_{\Gamma}^{\mu} M$ for
$\bar{\mu}:= \mu + Q \in \gh^*/Q$ (see Lemma \ref{lmnew}, (ii)).
If $M$ and $\Gamma$ are as in Proposition \ref{simple} then ${\cal
E}(M)$ contains all simple modules in ${\mathcal B}^{\chi}$ as
subquotients.

Some of the properties of the twisted localization are described
in the following proposition:

\begin{proposition} \label{locind}Let $M$ be a $\Gamma$-injective
$\gg$-module in $\mathcal B$.

$(i)$ ${\cal D}_{\Gamma} M \simeq M$ iff $M$ is
$\Gamma$-bijective.

$(ii)$ $\supp {\cal D}_{\Gamma}^{\mu} M = \mu+ \supp M +
\Span_{\Z} \Gamma$. Moreover, if $\nu_0 \in \supp M$ then $\dim
({\cal D}_{\Gamma}^{\mu} M)^{\nu'} = \max\, \{ \dim M^{\nu} \; |
\; \nu \in \nu_0 + \Span_{\Z}\Gamma\}$, whenever $\nu' \in \mu +
\nu_0 + \Span_{\Z} \Gamma$.

$(iii)$ Let $M$ be a module in ${\cal B}$ which has a unique
simple submodule. Then
  ${\cal D}_{\Gamma}^{\mu} M$ is
indecomposable whenever $M$ is indecomposable.
\end{proposition}
\begin{proof} Statement (i) is straightforward. (ii) follows from
a generalization of Lemma 4.4 in \cite{M}. Since $
\Phi^{\mu}_{\Gamma} \circ \Phi^{-\mu}_{\Gamma} = \Id $, to prove
(iii) is enough to show that ${\cal D}_{\Gamma} M$ is
indecomposable. Suppose ${\cal D}_{\Gamma} M = D_1 \oplus D_2$.
Then by our assumption $M$ has trivial intersection with one of
the modules $D_1$ or $D_2$, say $M \cap D_1 = 0$. We next show
that $(D_1)^{\nu'} = 0$, for a fixed $\nu' \in \supp {\cal
D}_{\Gamma} M$ (and thus $D_1=0$). We choose $\nu_0 \in \nu' +
\Span_{\Z}\Gamma$ such that $\dim M^{\nu_0} = \max\, \{ \dim
M^{\nu} \; | \; \nu \in \nu' + \Span_{\Z}\Gamma\}$. Then by (ii),
$M^{\nu_0} = ({\cal D}_{\Gamma} M)^{\nu_0} = (D_1)^{\nu_0} \oplus
(D_2)^{\nu_0}$ and therefore $(D_1)^{\nu_0} = 0$. However, (i)
implies that $D_1$ is $\Gamma$-bijective as a submodule of ${\cal
D}_{\Gamma} M$ and thus $(D_1)^{\nu'} = 0$.  \end{proof}

\begin{remark}
Statement (iii) of Proposition \ref{locind} remains valid if we
replace the condition $M \in {\mathcal B}$ by the weaker
requirement that $M$ is bounded in the $\Gamma$-directions only,
i.e. that the set $\{ \dim M^{\lambda}\; | \; \lambda \in
\lambda_0 + \Span_{\C}\Gamma\}$ is uniformly bounded for every
$\lambda_0 \in \supp M$.
\end{remark}

\begin{proposition} \label{mult}
Let $\gg = \gs \gp (2n)$ and $n>1$. Let $M$ be an indecomposable
$\gg$-module with unique simple submodule. Then there is a  set
$\Gamma$ consisted of $n$ commuting (i.e. orthogonal) long roots
such that $M$ is $\Gamma$-injective. Moreover, any composition
series of $M$ is multiplicity free, i.e. every two distinct simple
subquotients of $M$ are nonisomorphic.
\end{proposition}
\begin{proof} Let us prove the first statement. Suppose that there is a
long root $\beta$ for which both $f_{\beta}$ and $f_{-\beta}$ do
not act injectively on $M$. Let
$$
M_0:= \{m \in M\; | \; f_{\alpha}^N m = 0, \mbox{ some } N\}
\oplus \{m \in M\; | \; f_{-\alpha}^K m = 0, \mbox{ some } K\}.
$$
The sum is direct since for every simple $\gg$-module $P$, for
every $p \in P$, and for every long root $\beta$, we have that
$f_{\beta}^N p = f_{-\beta}^M p = 0$ implies $p=0$. The submodule
$M_0$ of $M$ is a direct sum of two nonzero submodules which
contradicts the initial assumption.

To prove the second statement choose $\mu \in \gh^*$ and $\Gamma
\subset \Delta$ so that $C:={\cal D}_{\Gamma}^{\mu}M$ is a
cuspidal module, i.e. all elements of $\gg \setminus \gh$ act
bijectively on $C$. We have that $C$ is semisimple (Theorem 1 in
\cite{BKLM}) and indecomposable (Proposition \ref{locind}, (ii)),
and hence it is simple. Let $N$ be the simple submodule of $M$.
Then $C\simeq {\cal D}^{\mu}_{\Gamma}N$, and therefore
$$
{\cal D}_{\Gamma}N \simeq \Phi^{-\mu}_{\Gamma}C \simeq
\Phi^{-\mu}_{\Gamma}{\cal D}_{\Gamma}^{\mu}M \simeq {\cal
D}_{\Gamma}M.
$$
By Proposition \ref{simple}, ${\cal D}_{\Gamma}M$ has a
multiplicity free compositions series, and so does its submodule
$M$.
\end{proof}

\begin{proposition} \label{projloc}
Let $\gg= \gs \gp (2n)$, $n>1$. Let $M$ be a simple module in
$\mathcal B$ and $\Gamma$ be a set of $n$ long roots such that $M$
is $\Gamma$-injective. Then ${\cal D}_{\Gamma}M$ and its
restricted dual $({\cal D}_{\Gamma}M)^*$ are the injective hull
and the projective cover of $M$ in ${\mathcal B}$, respectively.
\end{proposition}

\begin{proof} We first show that ${\cal D}_{\Gamma}M$ is
  injective, i.e. any exact sequence
$$
0 \to {\cal D}_{\Gamma}M \to M' \to N \to 0
$$
splits in $\mathcal B$. It suffices to prove this in the case when
$N$ is simple. Assume that a sequence does not split. Then $M'$
satisfies Proposition \ref{mult}. Since $\supp N \subset
\supp{\cal D}_{\Gamma}M$ and and $N$ has the same central
character as $M$, by Proposition \ref{simple}, $N$ is isomorphic
to some simple subquotient of ${\cal D}_{\Gamma}M$. Therefore $N$
is a subquotient of $M'$ with multiplicity higher than one, which
contradicts to Proposition \ref{mult}. The second statement
follows by duality. \end{proof}

\begin{corollary}\label{nice}Let $\gg= \gs \gp (2n)$, $n>1$. Then
  every simple object in $\mathcal B$ has a unique projective
  indecomposable cover and a unique injective hull.
\end{corollary}

\begin{remark} \label{sl2} Propositions \ref{mult} and
  \ref{projloc} are false for $n=1$. In fact, in this case the category
  $\mathcal B$ does not have injective and projective modules.
 To see this, let $\Omega$ denote the Casimir operator of $\gs \gl (2)$
 and $H$ be the standard element in the Cartan subalgebra. Let $P$ be an indecomposable
  projective module in $\mathcal B$, $M$ be some simple quotient and
$\mu \in \supp M$. There exists an integer $p$ and $\nu \in \mathbb C$
  such that $(\Omega -\nu)^p$ acts by zero on $P$.

 Let for $s \in \Z$, $I_s$ be the left ideal in  $U(\gg)$ generated by $H-\mu$ and
 $(\Omega -\nu)^s$, and let
  $F(s,\mu,\nu):=U(\gg)/I_s$. Then $\supp F(s,\mu,\nu)=\mu + Q$ and
  every weight has multiplicity $s$. Moreover,
  $F(s,\mu,\nu)$ is indecomposable with unique simple quotient
  isomorphic to $M$. Hence there exists a surjective homomorphism
$P \to F(s,\mu,\nu)$. However, if $s>p$ such homomorphism can not
be surjective which leads to a contradiction.
\end{remark}

\begin{example} \label{exproj}
\bigskip
Let $\mathcal {FIN}$ be the category of all weight $\gs \gp
(2n)$-modules with finite weight multiplicities and locally finite
action of the center of $U(\gg)$. It is not difficult to show that
every indecomposable module in $\mathcal{FIN}$ has finite length.
However, Corollary \ref{nice} does not hold if we replace the
category $\mathcal B$ by $\mathcal{FIN}$. Here is a
counterexample. Choose a parabolic subalgebra $\gp$ of $\gg$ such
that a Levi subalgebra $\gs$ of $\gp$ is isomorphic to $\gs \gl
(2)$. Choose $H \in \gs$, $\Omega \in U(\gs)$ and $\mu,\nu \in
\mathbb C$ as in the previous remark, so that $F(s,\mu,\nu)$ is a
simple $\gs$-module. Endow $F(s,\mu,\nu)$ with a structure of a
$\gp$-module by letting the radical to act by zero. Let
$$
M^s:=U(\gg) \otimes _{U(\gp)} F(s,\mu,\nu).
$$
 Then $M^s$ is indecomposable and belongs to $\mathcal {FIN}$. It is
 not difficult to see that
 $M^s$ has a unique simple quotient which we
 denote by $L$. We claim that $L$ does not have a projective cover in
 $\mathcal{FIN}$. This follows by reasoning similar to the one in
 the previous remark. Indeed, if $P$ is a projective cover of $L$,
 then there is a surjective map $P \to M^s$ for any $s$. Since $P$ has
 finite length, this is impossible.

\bigskip
\end{example}

\section{From bounded weight $ {\mathfrak s}{\mathfrak p}\left(2n\right)$-modules to weight
 $ {\protect \mathcal A}_{n} $-modules }

Let $ {\mathfrak g}={\mathfrak s}{\mathfrak p}\left(2n\right) $
with $ n\geq2 $. Every element $ X\in{\mathfrak g} $ can be
written in a block matrix form
\begin{equation} \left[
\begin{matrix}
A & B \\
C & -A^{t}
\end{matrix} \right]
\notag\end{equation} where $ A $ is an arbitrary $ n\times n
$-matrix, and $ B $ and $ C $ are symmetric $ n\times n
$-matrices. The maps
\begin{equation}
B \mapsto \sum_{i\leq j}b_{ij} t_{i}t_{j}\text{, }C \mapsto
\sum_{i\leq j}c_{ij} \partial_{i}\partial_{j} \notag
\end{equation}
can be extended to a homomorphism of Lie algebras
\begin{equation}
{\mathfrak g} \to {\mathcal A}_{n}
\notag\end{equation}
which induces a homomorphism
\begin{equation}
\omega\colon U\left({\mathfrak g}\right) \to {\mathcal A}_{n}.
\notag\end{equation} It is easy to see that the image of $ \omega
$ coincides with $ {\mathcal A}_{n}^{ev} $. If we fix the standard
basis of ${\mathcal A}_{n}$ we verify that $ \omega\colon
U\left({\mathfrak h}\right) \to {\mathcal A}_{n}^{0} $ is an
isomorphism. The representation of $ {\mathcal A}_{n}^{ev} $ in
the subspace $ W $ of even functions in $ {\mathbb
C}\left[t_{1},\dots ,t_{n}\right] $ is called the {\it Weil
representation}. One can check that $ W $ is irreducible. If $
I:=\operatorname{Ker} \omega $, then clearly $
I=\operatorname{Ann} W $ is a primitive ideal in $
U\left({\mathfrak g}\right) $. The center $ Z $ of $
U\left({\mathfrak g}\right) $ acts on $ W $ via the central
character $\sigma$ of $W$.

As follows from \cite{M}, for every simple module $ M $ in the
category $ {\mathcal B}^{\sigma} $,
\begin{equation}
\operatorname{Ann} M=I.
\notag\end{equation}

The next theorem follows from Proposition 12.1 in \cite{M}.

\begin{theorem} \label{translation} Let $\chi$ be a central character
  such that $\mathcal B ^{\chi}$ is non-empty. Then $\mathcal B ^{\chi}$ is
  equivalent to $\mathcal B ^{\sigma}$, with equivalence given
  by a translation functor.
\end{theorem}

(For the definition and properties of the translation functor see
\cite{Ja}.)

\begin{theorem} \label{th3}\myLabel{th3}\relax
Let $ M $ be any module from the category $ {\mathcal B}^{\sigma}
$. Then $ \operatorname{Ann} M=I $.

\end{theorem}

\begin{proof} It is sufficient to check that the statement holds for injective modules.
The latter follows form the fact that all injectives are obtained
via a localization as shown in Proposition \ref{projloc}.
\end{proof}

\begin{corollary} \label{cor3}\myLabel{cor3}\relax
The categories $ {\mathcal F}_{n}^{ev} $ and $ {\mathcal B}^{\chi}
$ are equivalent.

\end{corollary}

\section{Explicit description of all bounded $\gs \gp
(4)$-modules}\label{explic}

In this section we  explicitly describe all indecomposable objects
in $\cal B$ for $\gg := \gs \gp (4)$. We use the same notations as
in Section \ref{twloc}.

Let  $\Delta = \{ \pm \alpha_i, \pm \beta_i\, | \, i =1,2\}$ be
the root system of $\gg$ where  $\alpha_1, \alpha_2$, and
$\beta_1, \beta_2$ are the positive short and long roots,
respectively. Denote by $B: = \{ \alpha_1, \beta_2\}$  the
standard basis of $\Delta$ and let $\Gamma:=\{ \beta_1,
\beta_2\}$. There is an orthonormal basis
$\{\varepsilon_1,\varepsilon_2 \}$ of $\gh^*$ for which $\alpha_1
= \varepsilon_1 - \varepsilon_2$ and $\beta_2 = 2 \varepsilon_2$.
Let $W$ be the Weyl group of $\gg$, and let $s_{\alpha} \in W$
denote the reflection corresponding to the root $\alpha$.

For a $\gg$-module $M$ we denote by $M^*$ the restricted dual of
$M$. Note that $M^*$ is isomorphic to the twist
$M^{s_{\beta_1}s_{\beta_2}}$ of $M$ by $s_{\beta_1}s_{\beta_2} \in
W$. For a basis $B'$ of $\Delta$ and a weight $\lambda \in \gh^*$,
by $L_{B'}(\lambda)$ we denote the simple highest weight module
with highest weight $\lambda$ relative to the Borel subalgerbra
corresponding to $B'$. Put $\rho_{B'}$ for the half sum of the
$B'$-positive roots in $\Delta$.

For  $\gg$-submodules $A_1$ and $A_2$ of a $\gg$-module $A$, as
usual, the $A$-{\it diagonal} in $A_1 \oplus A_2$ is:
$$
D(A):= \{ (a, a) \in A_1 \oplus A_2 \; | \; a \in A\}.
$$
For the purpose of our construction we need a more general notion.
If $L$ is an endomorphism of $\C^k$, we define the $(A,L)$-{\it
diagonal} in $A_1^{\oplus k} \oplus A_2^{\oplus k}$ by
$$
D_L(A):= \{ (a, L(a)) \; | \; a \in A^{\oplus k}\}.
$$
In particular, for $k=1$ and $L = \Id$ we have $D_L(A) = D(A)$.

For $\eta \in \gh^*$ we set
$$
{\cal B}^{\chi}[\eta]:=\{M \in {\cal B}^{\chi}\; | \; \supp M
\subset \eta + Q\}.
$$
We next describe the simple objects of the subcategory ${\cal
B}^{\chi}[\eta]$ of ${\cal B}^{\chi}$. There are three types of
categories ${\cal B}^{\chi}[\eta]$ depending on the image $\eta +
Q$ of $\eta$ in the torus $\gh^*/Q$.

\medskip
$\bullet$ {\it Highest weight type:} $\eta +Q \in {\mathcal
H}{\mathcal W}(\chi)$. The simple objects of ${\cal
B}^{\chi}[\eta]$ are highest weight modules. There are two
elements $\eta +Q$ in $\gh^*/Q$ with this property. If
$L_B(\lambda^+)$ and $L_B(\lambda^-)$ are the two $B$-highest
weight modules in ${\cal B}^{\chi}$ then $\lambda^{\pm} + \rho_B =
m_1\varepsilon_1 \pm m_2\varepsilon_2$ for $m_i \in \frac{1}{2} +
\Z$ (note that $\lambda^- = s_{\beta_2}\lambda^+$). We fix
$\lambda^{\pm}$ so that $m_2 \geq -1/2$. Then the four highest
weight modules in ${\cal B}^{\chi}[\lambda^{\pm}]$ are:
$$
\begin{array}{c}
N^{\pm} :=
L_{s_{\beta_1}s_{\beta_2}(B)}(s_{\beta_1}s_{\beta_2}(\lambda^{\pm}))\\
W^{\pm} := L_{s_{\beta_2}(B)}(\lambda^{\pm} + \beta_2), \; \; \;
E^{\pm} :=
L_{s_{\beta_1}(B)}(s_{\beta_1}s_{\beta_2}(\lambda^{\pm}) -
\beta_2)\\
S^{\pm} := L_B(\lambda^{\pm})
\end{array}
$$
(standing for north, west, east, and south, respectively). Let
${\cal A^{\pm}}:=\{ N^{\pm}, E^{\pm}, S^{\pm}, W^{\pm}\}$. In
future we will consider modules either in ${\cal A^{+}}$ or in
${\cal A^{-}}$. For simplicity we will omit the superscripts and
will write ${\mathcal A}, N, W, E, S$.

\medskip
$\bullet$ {\it Cuspidal type:} $\eta +Q \in {\mathcal C}{\mathcal
U}{\mathcal S}{\mathcal P}(\chi)$. In this case there is only one
simple object in ${\cal B}^{\chi}[\eta]$ isomorphic to ${\cal
D}_{-\beta_1, -\beta_2}^{\eta - \lambda^+}L_{B}(\lambda^+)$, where
$\eta - \lambda^+ = x_1 \alpha_1 + x_2 \alpha_2$ with $x_i \notin
\Z$.

\medskip
$\bullet$ {\it Semi-plane type:} $\eta +Q \in {\mathcal
S}{\mathcal E}{\mathcal M}{\mathcal I}(\chi)$. There are two
simple objects in ${\cal B}^{\chi}[\eta]$ whose supports are
semi-planes. In this case $\eta +Q$ equals  $\lambda^+ + x
\varepsilon_1 +Q$ or $\lambda^+ + x \varepsilon_2 + Q$ for $x
\notin \Z $ which we will call {\it NW-ES type} and {\it NE-SW
type}, respectively. The two simple objects are isomorphic to
${\cal D}_{-\beta_1}^{\eta-\lambda^+}L_{B}(\lambda^+)$ and its
dual for the NW-ES type and ${\cal D}_{-\beta_2}^{\eta-\lambda^+}
L_{B}(\lambda^+)$ and its dual for the NE-SW type. Here $\eta -
\lambda \in  x \varepsilon_i + Q$ for $x \notin \Z$ and $i = 1$
(respectively, $i =2$) for the NW-ES (resp., NE-SW) type.
\medskip

\begin{example} In the special case when $\chi$ equals the central character
$\chi_0$ of the Weyl modules $L_B(\omega^+)$ or $L_B(\omega^-)$,
where $\omega^+ = - \frac{1}{2} \varepsilon_1 - \frac{1}{2}
\varepsilon_2$ and $\omega^- = - \frac{1}{2} \varepsilon_1 -
\frac{3}{2} \varepsilon_2$, we have that all simple objects in
${\cal B}^{\chi_0}$ have one-dimensional weight spaces. We may
simplify our considerations if we first restrict our attention to
the category ${\mathcal B^{\chi_0}}$ and then apply the
translation functor $\theta_{\chi_0}^{\chi} : {\mathcal
B^{\chi_0}} \to {\mathcal B^{\chi}}$, $\theta_{\chi_0}^{\chi}(M):=
\mbox{pr}_{\chi}(M \otimes L_B(\lambda^+ - \omega^+))$, where
$\mbox{pr}_{\chi}$ is the projection onto ${\mathcal B^{\chi}}$
(note that $L_B(\lambda^+ - \omega^+)$ is a finite dimensional
module). The highest weight part ${\mathcal
B^{\chi_0}}[\omega^{-}]$ of ${\mathcal B^{\chi_0}}$ is described
on Figure \ref{fig:1}. The other highest weight part, ${\mathcal
B^{\chi_0}}[\omega^{+}]$, can be pictured by rotating Figure
\ref{fig:1}  by $90^{\circ}$.

\end{example}

\begin{figure}[htbp]
  \begin{center}
    \leavevmode
  \setlength{\unitlength}{0.25in}

\begin{picture}(21,19)

\thicklines
\multiput(2.5,1.5)(0,1.5){11}{\circle{.25} }
\multiput(4,1.5)(0,1.5){11}{\circle{.25} }
\multiput(5.5,1.5)(0,1.5){11}{\circle{.25} }
\multiput(7,1.5)(0,1.5){11}{\circle{.25} }

\multiput(8.5,1.5)(0,1.5){11}{\circle{.25} }
\multiput(10,1.5)(0,1.5){11}{\circle{.25} }
\multiput(11.5,1.5)(0,1.5){11}{\circle{.25} }
\multiput(13,1.5)(0,1.5){11}{\circle{.25} }
\multiput(14.5,1.5)(0,1.5){11}{\circle{.25} }

\multiput(16,1.5)(0,1.5){11}{\circle{.25} }

\thicklines

\put(9.86,7.34){$\bullet$}
\put(9.5,6.8){$\omega^{-}$}

\put(9,0){$S^-$}

\put(9,17.5){$N^-$}

\put(17,9){$E^-$}

\put(1  ,9){$W^-$}

\put(16,16.5){\line(-1,-1){5.85}}

\put(16,15){\line(-1,-1){5.85}}

\put(10,10.5){\line(-1,0){1.5}}

\put(10,9){\line(1,-1){5.85}}

\put(8.5,9){\line(-1,-1){5.85}}

\put(10,7.5){\line(-1,0){1.5}}

\put(8.5,10.5){\line(-1,1){5.85}}

\put(8.5,9){\line(-1,1){5.85}}

\put(8.5,7.5){\line(-1,-1){5.85}}

\put(10,7.5){\line(1,-1){5.85}}

\end{picture}

\vspace{3ex}
    \caption{}
    \label{fig:1}
  \end{center}
\end{figure}

\begin{remark}
Let  ${\mathcal M}_{\chi}$ be the unique semisimple coherent
family with central character $\chi$ which is irreducible, i.e.
for which ${\mathcal M}_{\chi}[\lambda] := \oplus_{\mu \in \lambda
+ Q}({\mathcal M}_{\chi})^{\mu}$ is irreducible for some
$\lambda$. Another way to describe the three types of cosets $\eta
+ Q$ is via the generalized Shapovalov map $S_{\chi} : \gh^* \to
\C$ defined by $ \lambda \mapsto \det
(f_{\beta_1}f_{\beta_2}e_{\beta_1}e_{\beta_2})_{|{(\mathcal
M}_{\chi})^{\lambda}}$. We have that for $\eta \in \gh^*$ the zero
set of the restriction $S_{{\chi}|\eta +Q}$ of $S_{\chi}$ is
either empty, a line, or a union of two lines. These three cases
for $\eta +Q$ correspond to cuspidal, semi-plane, and highest
weight type, respectively.
\end{remark}

\medskip
\begin{lemma} \label{sp4loc}(i) The module ${\cal D}_{\beta_1}N$ (respectively, $ {\cal
D}_{\beta_2}N$) has length two and $({\cal D}_{\beta_1}N)/N \simeq
W$, (resp.,  $({\cal D}_{\beta_2}N)/N \simeq E$).

(ii) The module ${\cal D}_{\beta_1, \beta_2}N$ has length $3$ and:
$$
0 \subset L_1 \oplus L_2 \subset L_3 = ({\cal D}_{\beta_1,
\beta_2}N)/N
$$
where $L_1 \simeq E$, $L_2 \simeq W$, and $L_3/(L_1 \oplus L_2)
\simeq S$.
\end{lemma}

Let $T = (T_1, ..., T_k)$  be an ordered $k$-tuple of elements in
${\cal A}$. We call $T$ {\it admissible} if $T_i$ and $T_{i+1}$
are {\it successive} in ${\cal A}$, i.e. for $T_i=N$, we have
either $T_{i+1}=E$ or $T_{i+1}=W$, etc. For $X \in {\cal A}$ and
$T=(T_1, ..., T_k)$ for which $(X, T)$ is admissible we construct
 an indecomposable extension $X_T$ of $X$ for which $(X_T/X)^*
\simeq (T_1)_{T_2,...,T_k}$. A convenient way to represent $X_T$
is by a graph with a set of vertices $T\cup \{X\}$ and oriented
edges $T_{2i+1} \to T_{2i}$ and $T_{2i+1} \to T_{2i+2}$, $i \geq
0$, where $T_0:=X$. As an immediate application of Lemma
\ref{sp4loc} we define $W_{S,E} = E_{S,W} := {\cal D}_{\beta_1,
\beta_2}N$. In a similar way we set $N_{W,S}:= ({\cal D}_{\beta_1,
- \beta_2} E)/E, N_{E,S}:= ({\cal D}_{- \beta_1, \beta_2} W)/W$.

Since $W$ is a submodule of both $W_N$ and $E_{S,W}$ and $E$ is a
submodule of both $E_{S,W}$ and $E_N$ we  may define:
$$
N_{W,S,E}:= ((W_{S,E} \oplus W_N) / D_{i_W, j_W}(W))^*,
N_{E,S,W}:= ((E_{S,W} \oplus E_N) / D_{i_E, j_E}(E))^*.
$$
We might think of $ N_{W,S,E}$ as the $\beta_1$-localization of
the ``$W$-part'' of $W_{S,E}$. With similar reasoning we set:
$$
N^1 = N_{W,S,E,N}:= ((N_{W,S,E}^* \oplus E_N)/ D(E))^* \simeq
((N_{E,S,W}^* \oplus W_N)/ D(W))^*.
$$
We easily generalize the above constructions and for $X$ and $Y$
in  ${\cal A}$ define $X_{(Y,T)}$ using a "partial localization"
of $Y_T$. Also, if $T=(T_0, T_1)$ where $T_0$ has $l$ copies of
each element of ${\cal A}$ we set for simplicity $X_{T_1}^l:=X_T$
(we allow $T_1 = \varnothing$ as a $0$-tuple writing simply $X^l$
in this case). We  put also $X^0_0:= X$, $X^0_T:= X_T$ for a
$k$-tuple $T$, $0 \leq k \leq 3$,

We next notice that $E$ and $W$ are submodules of $W_{N, E}$ and
$W_{S,E}$, so for every $c \in \C$ and a positive integer $k$ we
define
$$
N^k_{\lambda}:=(W_{N, E}^{\oplus k} \oplus W_{S, E}^{\oplus
k})/(D_{\Id}(W^{\oplus k}) \oplus D_{J^k_{c}}(W^{\oplus k})),
$$
where $J^k_{c} \in End(\C^k)$ is represented by a single Jordan
block with $c$ on the diagonal. Note that $N_0^k \simeq
N_{E,S,W}^{k-1}$ for $k \geq 1$.

In similar fashion we construct $X_{c}^k$ for every $X$ in  ${\cal
A}$. We set $A_{c}^l :=N_{c}^l \simeq S_{c}^l $ and $B_{c}^l
:=E_{c}^l \simeq W_{c}^l $, $l \geq 1$. Finally, denote by $P_X$
the projective cover of $X$. Note that $P_X$ is the
$\Gamma_X$-localization of $X$ where $\Gamma_X$ is the set of
those two long roots for which $X$ is $\Gamma_X$-localizable.

\begin{proposition} Up to an isomorphism,
the complete list of the indecomposable objects in ${\cal
B}^{\chi}$ includes:

(i) Highest weight type: $P_X$, $X_T^k$, $(X_T^k)^*$, $A_{c}^l$,
$B_{c}^l \simeq (A_{c}^l)^*$, where $X \in {\cal A}$, $T$ is an
$n$-tuple, $0 \leq n \leq 3$, $k \geq 0$, $l\geq 1$ $c \in \C$. Up
to a twist of an element of the Weyl group we have five types (the
projective and four series) of modules: $P_N$, $N^k, N_E^k,
N_{E,S}^k,$ and $N_{c}^k$.

(ii) Cuspidal type: ${\cal D}_{\beta_1, \beta_2}^{\mu}N$ with $\mu
= x_1 \alpha_1 + x_2 \alpha_2$, $x_i \notin \Z$.

(iii) Semi-plane type: ${\cal D}_{\beta_1}^{\mu}N$, $({\cal
D}_{\beta_1}^{\mu}N)^*$, ${\cal D}_{\beta_2}^{\nu}N$, $({\cal
D}_{\beta_2}^{\nu}N)^*$, ${\cal D}_{\beta_1, \beta_2}^{\mu}N$,
$({\cal D}_{\beta_1, \beta_2}^{\mu}N)^*$, ${\cal D}_{\beta_1,
\beta_2}^{\nu}N$, $({\cal D}_{\beta_1, \beta_2}^{\nu}N)^*$, where
$\mu = x_1 \varepsilon_1 + x_2 \varepsilon_2$ and $\nu = y_1
\varepsilon_1 + y_2 \varepsilon_2$ are such that $x_1 \notin \Z$,
$x_2 \in \Z$, $y_1 \in \Z$, $y_2 \notin \Z$. Up to a twist of the
Weyl group there are two types: ${\cal D}_{\beta_1}^{\mu}N$ and
${\cal D}_{\beta_1, \beta_2}^{\mu}N$.
\end{proposition}

\def\cprime{$'$} \def\cprime{$'$} \def\cprime{$'$} \def\cprime{$'$}
  \def\cprime{$'$}
\providecommand{\bysame}{\leavevmode\hbox
to3em{\hrulefill}\thinspace}
\providecommand{\MR}{\relax\ifhmode\unskip\space\fi MR }
\providecommand{\MRhref}[2]{%
  \href{http://www.ams.org/mathscinet-getitem?mr=#1}{#2}
} \providecommand{\href}[2]{#2}

\end{document}